\documentclass{etds}

\ETDS{31}{6}{1865–1887}{2011}

\newtheorem{theo}{Theorem}
\newtheorem{prop}{Proposition}[section]
\newtheorem{lemma}{Lemma}[section]

\begin{document}
\runningheads{T.\ Yarmola}{A pathological random perturbation of the Cat Map}

\title{An example of a pathological random perturbation of the Cat Map}

\author{TATIANA YARMOLA}

\address{Department of Mathematics, Mathematics Building, University of Maryland,
College Park, MD 20742-4015, USA \\
\email{yarmola@math.umd.edu}}

\recd{$6$ April $2010$}

\begin{abstract}
In this paper we give an example of a random perturbation of the Cat Map that produces a \textquotedblleft global statistical attractor" in the form of a line segment. The transition probabilities for this random perturbation are smooth in some but not all directions. All initial distributions on $\mathbb{T}^2$ are attracted to distributions supported on this line segment.
\end{abstract}

\section*{Introduction.}

Random perturbations of dynamical systems are important in modeling noise and other types of uncontrolled fluctuations. Given a Riemannian manifold $M$ and a mapping $f:M \to M$, a random perturbation of $f$ is defined as a Markov chain on $M$ such that for every $x$, the transition probability $P(\cdot|x)$ is given by $Q_{fx}$ where $Q_{fx}$ is a probability distribution
that depends only on $f(x)$ and is not far from the point mass at
$f(x)$. Intuitively it means that a particle jumps from $x$ to $f(x)$ and then disperses randomly near $f(x)$ with the distribution $Q_{fx}$.

Assuming $M$ is compact and $Q_{fx}$ depend continuously on $x$, such Markov chain admits stationary measures. We will refer to these stationary measures as invariant measures. If we denote the Riemannian measure on $M$ by $m$, an important question to ask is whether the system admits invariant measures absolutely continuous with respect to $m$.
Regarding such invariant measures as \textquotedblleft natural" invariant measures (to distinguish them from any singular measures that may exist), we may ask if such measures are unique and if so, whether they capture the asymptotic dynamics of almost all points.
The answer to these questions is straightforward if we assume that $\{ Q_{fx} \}$ are absolutely continuous with respect to $m$  (i.e. have a density) for every $x$. Given the absolute continuity of $\{ Q_{fx} \}$, every invariant
measure of the perturbed dynamics is also absolutely continuous with respect to $m$.

Interesting problems arise for degenerate random perturbations in which the probability
distributions $Q_{fx}$ do not necessarily have densities. This subject has not yet been carefully studied, although there are many important applications. In many real systems, perturbations do not occur everywhere or uniformly in all directions. Frequently, irregular patches on the domain may introduce random patterns or perturbations can occur on the boundary. Common examples that can admit such perturbations include billiards, hard ball systems and chains of Hamiltonian systems coupled to heat baths.

For definiteness, let us consider uniform rank one perturbations, i.e. assume that the probability distributions $Q_x$ are uniform and supported on $1$-dimensional disks of length $2\epsilon$ centered at $x$. This assumption is not essential for many general properties of the perturbed system, but allows to describe things more precisely. Since we have to require $Q_x$ to depend continuously on $x$ in order to ensure the existence of the invariant measures, it makes sense to assume that the supporting intervals lie along the flow defined by some $C^\infty$ vector field.

The following argument suggests that, provided the dynamics are rich enough, it is natural to expect that there should exist an absolutely continuous invariant measure for the system subject to a degenerate random perturbation. If we start with the point measure at $x$, $\delta_x$, and push it forward by the perturbed dynamics, after the first step the new measure is equal to $Q_{fx}$ and is supported on a $1$-dimensional curve. If the image of
this curve under $f$ is not tangent to the vector field at any point, then the
\textquotedblleft smearing" will produce a  transition probability
$P^2(\cdot |x)$ that has a two-dimensional density, i.e. it acquires
an extra dimension from the perturbation. Now consider the $f$-image
of the support of $P^2(\cdot |x)$. If this image is never tangent to the
vector field, then in the next step yet another dimension is
acquired, i.e. $P^3(\cdot |x)$ now has a $3$-dimensional density. This
process may be continued as long as the non-tangency condition is
satisfied. Thus starting with the point measure $\delta_x$, we reach an $n$-dimensional density in $n$ steps. It could be shown that if $f$ is a $C^2$ diffeomorphism and for every $x$, there exists $k$ such that $P^k(\cdot |x)$ is absolutely continuous with respect to $m$, then every invariant measure of the perturbed system has a
density \cite{AnosovPaper}.

On the other hand, if for some points $x$ the $P^n(\cdot |x)$ fail to
acquire a density, in many circumstances the perturbed system will have singular invariant measure(s). Such invariant measures may or may not coexist with absolutely continuous invariant measure(s). Their supports may be all of the phase space or they may be quite small. We say that a singular invariant measure is \emph{invisible} if except for a Lebesgue measure zero set of initial conditions, such an invariant measure cannot be reached asymptotically. We say that a singular invariant measures is \emph{visible} otherwise. We are not especially interested in singular measures that are invisible, that is to say, we do not regard the situation as pathological if for Lebesgue a.e. initial condition $x$, $P^n(\cdot|x)$ eventually converges to an absolutely continuous invariant measure even if singular measures may exist. Invariant measures that are visible can be viewed as being analogous to SRB measures for deterministic systems.

When a singular measure $\mu$ attracts a Lebesgue positive measure set of initial conditions we call its support a \emph{statistical attractor}. If for every $x\in M$, $\frac{1}{n} \sum_{i=0}^{n-1} P^i(\cdot|x) \to \mu$), we will call the support of $\mu$ a \emph{global statistical attractor}.

The following question motivated the example presented in this paper:

\emph{Suppose the perturbation is of rank one and the dynamics of f are rich in a sense of having strong hyperbolicity and mixing properties. Assume further that almost all points acquire density in the sense discussed above. Is it still possible for a statistical attractor or a global statistical attractor to occur? }

We answer this question in the affirmative by proving the existence of a global statistical attractor for a rank one perturbation of linear hyperbolic toral automorphism. For definiteness, we will work with the Cat Map and produce a concrete example of a rank one random perturbation with a global statistical attractor in the form of a line segment (which is a local stable manifold of the fixed point). The set of parameters for which our construction possesses a global statistical attractor is open except for the requirement that the vector field has to be parallel to the stable direction one some local stable manifold of the fixed point. Our construction can be generalized to other Anosov Diffeomorphisms.

\section{General Idea.}

\subsection{Settings.}

Let $f:\mathbb{T}^2 \to \mathbb{T}^2$ be a $C^2$-diffeomorphism, $V$ a unit $C^\infty$ vector field on $\mathbb{T}^2$ and $\epsilon$ some fixed number. By $m$ we will denote the Riemannian measure on $\mathbb{T}^2$.

For each $x \in \mathbb{T}^2$, define distributions $Q_x$ as follows:\\
Let $\gamma_x$ be a curve in $\mathbb{T}^2$ along the flow of the vector field $V$ such that $x\in \gamma_x$, and for any $y\in \gamma_x$, the Riemannian distance from $x\in \gamma_x$ to any $y\in \gamma_x$ is $\leq \epsilon$. Since $\mathbb{T}^2$ is a Riemannian manifold and $V$ is $C^\infty$, $\gamma_x$ is well defined. Let $Q_x$ be the uniform distribution on $\gamma_x$.

By the random perturbation of $f$, given $\{Q_x\}$, we will mean the Markov chain $X_n$, $n=0,1,2, \cdots$ with transition probabilities $P(A|x)=P\{X_{n+1} \in A: X_n=x\} = Q_{fx}(A)$ defined for any $x\in \mathbb{T}^2$ and Borel set $A \subset \mathbb{T}^2$ \cite{kifer}. Given $f:\mathbb{T}^2 \to \mathbb{T}^2$ and unambiguously defined family $\{Q_x\}$ we will denote the randomly perturbed dynamics with $P(A|x)=Q_{fx}$ by $\mathcal{F}$. This notation is going to be used for short reference of the dynamics and to indicate pushing measures forward by the dynamics in a sense that $\mathcal{F}_*\nu = \int_{\mathbb{T}^2}P(A|x)d\nu$.

\proc{Definition.}
A probability measure $\mu$ on $M$ is called an invariant measure of the Markov chain $X_n$ if for any Borel set $A \subset \mathbb{T}^2$, $\int_{\mathbb{T}^2} P(A|x) d\mu(x) = \mu(A)$. By the invariant measure of a random perturbation of $f$ we will mean the invariant probability measure of the Markov chain with the corresponding transition probabilities.
\medbreak

\proc{Remark.} By Prop 1.4 in \cite{kifer} our perturbed system has at least one invariant measure.
\medbreak

\subsection{Perturbations on $\mathbb{T}^2$}
Let $f$ be the Cat Map defined by iterations of $\left(
      \begin{array}{cc}
        2 & 1 \\
        1 & 1 \\
      \end{array}
    \right)$-matrix mod 1. Then $\mu=\frac{3-\sqrt{5}}{2}$ and $\lambda = \frac{3+\sqrt{5}}{2}$ are the contraction and the expansion rates along the stable and unstable manifolds.

\begin{figure}
  \centering{\scalebox{0.5}{\includegraphics{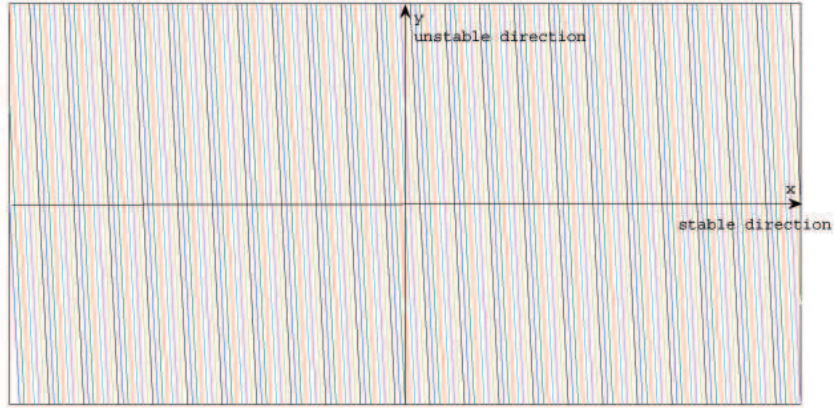}}}
  \caption{Integral lines of the \textquotedblleft tilted" vector field.}\label{fig:straight vector field}
\end{figure}

Consider the unit vector field $V_0$ on $\mathbb{T}^2$ everywhere parallel to the unstable direction. Although the Riemannian measure on $\mathbb{T}^2$ is clearly invariant for the perturbed dynamics $\mathcal{F}$, we do not have a pattern of \textquotedblleft acquiring density" in a sense it is described in the introduction and the fact that an invariant measure is absolutely continuous could be considered as a pure coincidence. If we \textquotedblleft tilt" $V_0$ slightly to form a small angle $\varphi$ with the unstable direction (see Figure \ref{fig:straight vector field}), the pattern of \textquotedblleft acquiring density" works everywhere and the Riemannian measure is still invariant for this situation. This is the starting point of our construction.

\begin{figure}
  \centering{\scalebox{0.5}{\includegraphics{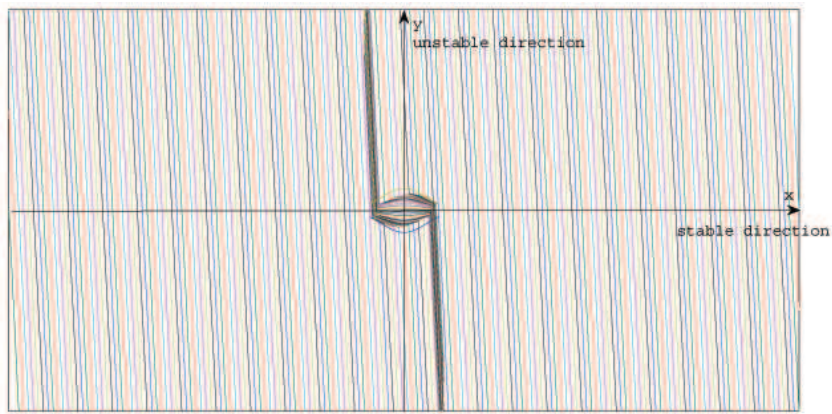}}}\\
  \caption{Integral lines of the \textquotedblleft kinked" vector field.}\label{fig:gaussian vector field}
\end{figure}

Our example is obtained by creating a \textquotedblleft kink" to this \textquotedblleft tilted" vector field in a small neighborhood $U'$ of the fixed point $(0,0)$. Specifically we will \textquotedblleft bend" the the vector field so that it is parallel to the stable direction on a short segment of the local stable manifold of the fixed point. See Figure \ref{fig:gaussian vector field}.

More precisely, set the coordinate frame with the stable manifold parallel to the $x$-axis and the unstable - to the $y$-axis with $(0,0)$ representing the fixed point. Then the dynamics of the Cat Map is described locally by $f(x,y)=(\mu x, \lambda y)$ around $(0,0)$. Define $W \subset U' \subset \mathbb{T}^2$ to be neighborhoods around the fixed point bounded by the Gaussian-shaped curves $y=\pm b\exp\{-\frac{x^2}{2\sigma^2}\}$, $y=\pm(b+\beta)\exp\{-\frac{x^2}{2\sigma^2}\}$ and line segments $x=a$, $x=a+\frac{\beta}{2}$ respectively for some $a,b, \beta, \sigma >0$. We are going to require that the new vector field is everywhere parallel to a family of Gaussian-Shaped curves $y=Ke^{-\frac{x^2}{2\sigma^2}}$ in $W$ and is $C^\infty$ extended in $U' \setminus W$ to match the \textquotedblleft tilted" vector field in $(U')^c$.

Clearly the random dynamics defined with the new vector field still has a pattern of \textquotedblleft acquiring density" everywhere except a short segment of the local stable manifold of the fixed point.

It turns out that under certain conditions on parameters, a perturbed vector field of this kind can break all the nice properties of the Cat Map dynamics and lead to the existence of a global statistical attractor, i.e. there exists a singular measure $\mu$ such that for every $x\in M$, $\frac{1}{n} \sum_{i=0}^{n-1} P^i(\cdot|x) \to \mu$. The global statistical attractor in this case, the support of $\mu$, is a short segment of the stable manifold around the fixed point. The set of parameters for which this happens is open except for the fact that the vector field must align with the stable direction on some interval around the fixed point.

\begin{theo}
Let $f$ be the Cat Map. Then for an open set of parameters defining the random dynamics $\mathcal{F}$

\begin{enumerate}
\item The only $\mathcal{F}$-invariant measure is the unique singular measure $\mu$ supported on the segment $[-a,a] \times \{0\}$ of the local stable manifold of the fixed point in the $E^s-E^u$-coordinate frame.
\item If $\nu$ is any Borel probability measure on $\mathbb{T}^2$, then $\mathcal{F}^n_* \nu$ has all its limit points supported on $[-a,a] \times \{0\}$. In particular, $\frac{1}{n}\sum_{i=0}^{n-1} \mathcal{F}_*^i \nu \to \mu$.
\end{enumerate}
\end{theo}

The mechanism that leads to this phenomena is as follows: without the random perturbation, orbits that come near the fixed point move away from it following the unstable direction. Our random perturbation as described above \textquotedblleft smears" them out roughly in the stable direction when $\epsilon$ is large enough. Since the integral curves of the vector field are Gaussian-shaped in $W$, part of the mass may also be pulled closer to $x$-axis than it was before. Our main result is that when the tendency to be pulled back to the $x$-axis (due to smearing) is stronger than the tendency to move away from it (due to the unperturbed Cat Map dynamics), a global statistical attractor will result.

In section 2 we will describe the local picture and show that it achieves the goal of creating a global statistical attractor under some artificial assumption about the return rates for the global dynamics. In section 3 we will show that this assumption is very similar to the dynamics of the Cat Map and show that with the appropriate modifications, the argument in section 2 carries through for the real perturbation of Cat Map.

\section{Local Analysis}

Let $W\subset \mathbb{T}^2$ be an open neighborhood of the fixed point enclosed by
the curves $x=\pm a$ and $y=\pm be^{-\frac{x^2}{2\sigma^2}}$ for
some $a,b, \sigma >0$. We consider in this section $f:\mathbb{T}^2\to \mathbb{T}^2$ defined as follows: for $(x,y)\in W$, $f(x,y)=(\mu x, \lambda y)$,
where $0<\mu=\frac{3-\sqrt{5}}{2}<1$, $\lambda=\frac{3+\sqrt{5}}{2}>1$ and $f$ is $C^\infty$ extended beyond W in
some fashion that we will describe later.

Consider the family of Gaussian-shaped curves, $\{\gamma_K\}$, given
by $\gamma_K(x) = Ke^{-\frac{x^2}{2\sigma^2}}$, where $\sigma$ is
the same as above and $K$ is the level of the curve determined by the
coordinate of $y$ when $x=0$. Since through any $p\in W$ passes a
unique curve of this type, we can consider a continuous unit vector
field $V$ tangent to $\{\gamma_K\}$ in $W$. Fixing $V(0,0)=(1,0)$ makes $V$ well defined when extended to the whole $\mathbb{T}^2$.

\begin{figure}[center]
  \centering{\scalebox{0.7}{\includegraphics{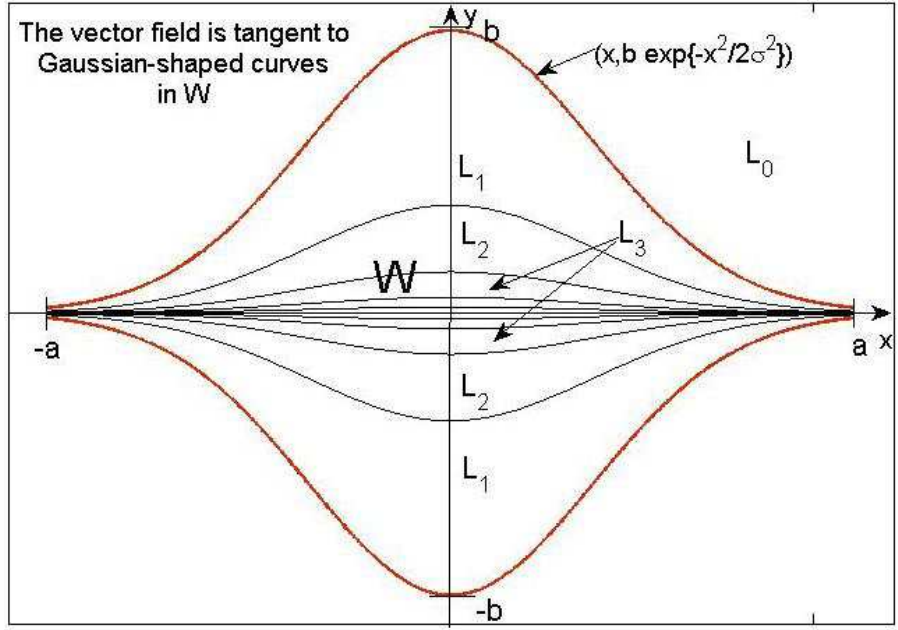}}}\\
  \caption{Neighborhood W.}\label{fig:W}
\end{figure}

We are going to add a $1$-dimensional random perturbation to $f$
prescribed by $V$ and some $\epsilon>0$. More
precisely, $\forall p\in W$, let K be such that $p \in \gamma_K$.
Define $I_\epsilon(p) = \{ s \in \gamma_K: d_K(p,s)<\epsilon\}$, where $d_K$
is the Riemannian distance along $\gamma_K$. Let $Q_p$ be the uniform
distribution on $I_\epsilon(p)$. Define $\mathcal{F}$ as in section 1.2 with $P(A|p)=Q_{fp}(A)$.

Let us tie $a$ and $\epsilon$ together to satisfy $a=\frac{\epsilon}{1-\mu}$. A simple computation ensures that all possible random images of $[-a,a] \times \{0\}$ belong to $[-a,a] \times \{0\}$ and therefore there exists an invariant measure supported on the interval $[-a,a] \times \{0\}$. Moreover, the push forward of any measure $\nu$ supported on $[-a,a] \times \{0\}$ is absolutely continuous with respect to the Lebesgue measure on $[-a,a] \times \{0\}$ with density values less or equal to $\frac{1}{2\epsilon}$ everywhere. Therefore any any invariant measure on $[-a,a] \times \{0\}$ is absolutely continuous with respect to the Lebesgue measure on $[-a,a] \times \{0\}$. It is also easy to see the invariant measure invariant measure on $[-a,a] \times \{0\}$ must be unique.

We are going to show that, under certain assumptions on $f$ and $V$ on $\mathbb{T}^2 \setminus W$, $\epsilon$, $a$, $b$, and $\sigma$, all the limit points of $\mathcal{F}^n_*\nu$ are singular (not necessarily invariant) measures supported on $W \cap \{ x-axis \}$. In addition, $\frac{1}{n}\sum_{i=0}^{n-1} \mathcal{F}^i_*\nu$ converges weakly to the unique singular invariant measure supported on $W \cap \{ x-axis \}$ (i.e. interval $[-a,a]$ of $x$-axis).

\subsection{Assumptions.}

Suppose $f|_W$ and $V|_W$ are defined as above. In order to create a statistical attractor we will need two provisional assumptions:
\begin{itemize}

\item (A1) Suppose $\epsilon$ is large enough, precisely
$$\epsilon > 3(\sigma \sqrt{\frac{4ln(\lambda)}{1-\mu^2}} + b) $$

\item (A2) $f|_W$ and $V|_W$ are
$C^\infty$ extended outside $W$ in some fashion to satisfy the
following:\\
\emph{For any finite measure $\nu$ on $\mathbb{T}^2$,
$(\mathcal{F}_*(\nu|_{W^c}))(W) \geq q\nu(W^c)$, where $q>0$ is some constant
and $\nu|_{W^c}(A)=\nu(A\cap W^c)$. In other words,
we assume that for any measure $\nu$, $q$-fraction of $\nu|_{W^c}$
gets to $W$ each step the measure is pushed forward by the perturbed
dynamics $\mathcal{F}$.}

\end{itemize}

\proc{Remark.} The assumption (A2) is not realistic in general because it requires a fraction of the measure outside $W$ to return to $W$ each time the measure is pushed forward by the perturbed dynamics. On the other hand, a condition along these lines is satisfied by the Cat Map. See section 3.
\medbreak

\begin{prop} \label{prop:local}
If $f$, $V$ and $\epsilon$ are defined as above and satisfy (A1) and (A2), then the only $\mathcal{F}$-invariant measure is singular and
is supported on the interval $[-a,a]$ of the x-axis. Moreover, if $\nu$ is any Borel probability measure on $\mathbb{T}^2$, then $\mathcal{F}^n_* \nu$ has all its limit points supported on $[-a,a] \times \{0\}$.
\end{prop}

\proc{Idea of Proof.} We are going to decompose $\mathbb{T}^2 \setminus
([-a,a]\times \{0\})$ into certain layers $L_n$ and compare the perturbed
dynamics $\mathcal{F}$ on these layers with a transient countable state
Markov chain. That will lead to the conclusion that all the invariant
measures for $\mathcal{F}$ must be supported on $[-a,a]\times \{0\}$.

\subsection{Dynamics on the layers}.\\

\proc{Definition.} \label{def:layer}
Let $c_0=b$, $c_1=\frac{1}{\lambda}c_0$,
$c_n=\frac{1}{\lambda}c_{n-1}$.
Define the \emph{layers} $$L_n = \{(x,y):
|x|\leq a, c_n\exp\{-\frac{x^2}{2\sigma^2}\}\leq |y| <
c_{n-1}\exp\{-\frac{x^2}{2\sigma^2}\}\}$$ for $n \geq 1$ and $L_0=W^c$. See Figure \ref{fig:W}.
\medbreak

\proc{Definition.} \label{def:level}
Define the \emph{level} of a Gaussian-shaped curve $\gamma_K(x) =
K\exp\{-\frac{x^2}{2\sigma^2}\}$ to be its value at $x=0$, namely K.
For any $(x,y) \in W$, let K(x,y) be the
\emph{level} of the unique Gaussian-shaped curve $\gamma_K(x) =
K\exp\{-\frac{x^2}{2\sigma^2}\}$ passing through $(x,y)$.
\medbreak

We are going to see now how $\mathcal{F}$ induces the dynamics on the layers
$L_n$.

\begin{lemma} \label{lemma:dynamics on layers}
If $(x,y) \in L_n$, $n \geq 1$, then $I_\epsilon(f(x,y)) \subset \cup_{k=n-1}^{\infty} L_k$, i.e. after one application of $\mathcal{F}$ any random image of $(x,y)$ moves
\textquotedblleft away" from the $x$-axis by at most one layer.\\
Moreover, if $(x,y) \in L_n$, $n \geq 1$, and $x \not\in [-x_0,x_0]$, then $I_\epsilon(f(x,y)) \subset \cup_{k=n+1}^{\infty} L_k$, where $x_0=\sigma
\sqrt{\frac{4ln(\lambda)}{1-\mu^2}}$. I.e. after one application of $\mathcal{F}$ any random image of $(x,y)$ moves \textquotedblleft towards" the $x$-axis by at least one layer.
\end{lemma}

\proc{Proof.} Simple computations yields:
\begin{itemize}
\item If $(x,y)\in L_n$, $n \geq 1$, then $f(x,y)\in \cup_{k=n-1}^{\infty}L_k$

\item If $x \not\in [-\sigma
\sqrt{\frac{4ln(\lambda)}{1-\mu^2}},\sigma
\sqrt{\frac{4ln(\lambda)}{1-\mu^2}}] := [-x_0,x_0]$, then
$K(f(x,y))\leq \frac{1}{\lambda}K(x,y)$ $\Rightarrow$\\
if $(x,y)\in L_n$ and $x \not\in [-x_0,x_0]$, then $f(x,y)\in \cup_{k=n+1}^{\infty}L_k$.

\item If $(x,y)\in W$ and $f(x,y)\in L_k$, $k \geq 1$, then $I_\epsilon(x) \subset L_k$ since we perturb by $\epsilon$ along the Gaussian-shaped curves and the $x$-coordinate cannot get beyond
$\pm\frac{\epsilon}{1-\mu}$.
\end{itemize}

\ep \medbreak

Now we are ready to discuss how the dynamics pushes measures forward between the layers.

\proc{Definition.} \label{def:layer index}
Let $\mu$ be any finite measure in $W$. We are going to
say that at least $\alpha$-fraction of $\mu$ \emph{increases the layer index}
under $\mathcal{F}$ if
$$ \sum_{k=1}^\infty [\mathcal{F}_*(\mu|_{L_k})](\cup_{i \geq k+1} L_i) \geq \alpha \mu(W) $$.
\medbreak

The idea of producing a singular limit of the pushed forward measures
is to increase $\epsilon$ such that measure spreads along $x$-direction far
enough to make a big fraction of it \emph{increase the layer index} in the
subsequent steps. This effect can also be achieved by decreasing
$\sigma$, the \textquotedblleft standard deviation" parameter in the Gaussian-shaped
curves. That is an important advantage since in many specific
examples we cannot increase $a=\frac{\epsilon}{1-\mu}$ indefinitely,
while scaling by $\sigma$ will ensure that we still get the effect
we want.

\begin{lemma} \label{lemma:epsilon}
Suppose $\nu$ is distributed uniformly on a piece of the level $K$ Gaussian-shaped curve, $\gamma_K$, of length $2\epsilon$ and symmetric with
respect to $y$-axis, where $K < b$. If
$$\epsilon > 3(\sigma \sqrt{\frac{4ln(\lambda)}{1-\mu^2}} + b) $$
then at least $\frac{2}{3}$ of $\nu$ will \emph{increase the layer index} when
pushed forward by the perturbed dynamics.
\end{lemma}

\proc{Proof.}
Let us estimate the amount of the uniformly distributed measure on $\gamma_K$ within the bounds $[-x_0,x_0]$:

$$\int_0^{x_0} \sqrt{1+ (\frac{dK\exp\{-\frac{x^2}{2\sigma^2}\}}{dx})^2}dx = \int_0^{x_0} \sqrt{1+ \frac{x^2K^2}{\sigma^4}\exp\{-\frac{x^2}{\sigma^2}\}}dx \leq $$
$$\int_0^{x_0}(1+\frac{xK}{\sigma^2}\exp\{-\frac{x^2}{2\sigma^2}\})dx = x_0+K(1-\exp\{-\frac{x^2}{2\sigma^2}\}) \leq \sigma \sqrt{\frac{4ln(\lambda)}{1-\mu^2}} + K < \frac{\epsilon}{3}$$

Using Lemma \ref{lemma:dynamics on layers} we conclude that at least $\frac{2}{3}$-fraction of $\nu$ will \emph{increase the layer index} when
pushed forward by the perturbed dynamics. \ep \medbreak

\begin{lemma} \label{lemma:2/3}
Let $\nu$ be any finite measure supported on $W \cap f^{-1}W$,
$\epsilon$ as above. Then $\mathcal{F}_*\nu$ is supported in W
and is spread in the $x$-direction enough that for each further iteration  at
least $\frac{2}{3}$ of the measure $(\mathcal{F}^n_*\nu)|_W$, $n > 1$,
\emph{increases the layer index} when pushed forward by the perturbed
dynamics $\mathcal{F}$.
\end{lemma}

\proc{Proof.} Let $\varrho=(\mathcal{F}^{n}_*\nu)(W)$, $n>1$.
Because the perturbation occurs for every iteration of the
dynamics, the one dimensional conditional density of $(\mathcal{F}^{n}_*\nu)|_W$ along any Gaussian-shaped curve $\gamma_K$ cannot exceed the value of
$\frac{1}{2\epsilon}$ at any point. Thus
$$(\mathcal{F}^{n+1}_*\nu)(W \cap ([-x_0,x_0]\times [-b,b])) \leq \varrho \frac{1}{2\epsilon} \frac{2\epsilon}{3} = \frac{\varrho}{3}.$$
Therefore at least $\frac{2}{3}$-fraction of $(\mathcal{F}^n_*\nu)|_W$ \emph{increases the layer index} under $\mathcal{F}$. \ep \medbreak

\subsection{Dominating Markov chain}

In this section we are going to compare the dynamics induced on the layers $L_n$ with a transient countable state Markov chain. This would allow us to conclude that no invariant measures are supported on $\cup_{n=0}^{\infty}L_n$, which implies that the invariant measure supported on $[-a,a] \times \{0\}$ is unique. Markov chain transience would also lead to the conclusion that any initial distribution converges to a distribution supported on $[-a,a] \times \{0\}$.

Let $S_0, ... S_n,  ...$ be a countable state Markov chain with the
following transition probabilities: for $k \geq 1$,
$P(k,k+1)=\frac{2}{3}$, $P(k,k-1) = \frac{1}{3}$; $P(0,1)=q$ and
$P(0,0)=1-q$. For all other $k$, $l$,
let $P(k,l)=0$. This defines an irreducible birth and death chain.

\begin{lemma} \label{lemma:mc transient}
The Markov chain $S_0, ... S_n,  ...$ is transient.
\end{lemma}

The Markov chain $S$ has transition probabilities biased to move the dynamics towards the states with larger indices: for the states indexed with $k \geq 1$, the probability of increasing the index is twice the probability of decreasing it. Lemma \ref{lemma:mc transient} is an easy consequence of the following probability exercise:

\begin{lemma} \label{lemma:durett} \cite[section 5.3, Exercise 3.7]{Durrett} \\
Let $X_0, X_1, \cdots, X_n, \cdots$ be an irreducible countable state Markov chain,
$\varphi$ a nonnegative function satisfying the following:

\begin{itemize}
\item $\varphi(x) \to 0$ as $x \to \infty$.
\item $\exists$ finite nonempty set F such that
           \begin{itemize}
           \item $\varphi(x)>0$ for $x\in F$  and
           \item $E_x\varphi(X_1)\leq \varphi(x)$ for $x\not \in F$.
           \end{itemize}
\end{itemize}
Then the chain is transient.
\end{lemma}

\proc{Proof of Lemma \ref{lemma:mc transient}.}
Let $\varphi(k)=\frac{1}{2^k}$. Then $\varphi(k)>0$ everywhere and
$\varphi(k) \to 0$ as $k \to \infty$. The Markov chain $S_0, \cdots,
S_n, \cdots $ has the property that for $k \geq 1$,
$$\varphi(k) = \frac{1}{2^k} = \frac{1}{2^k}(\frac{1}{3}+\frac{2}{3}) =
\frac{2}{3}\frac{1}{2^{k+1}} + \frac{1}{3}\frac{1}{2^{k-1}}= E_k
\varphi(S_1)$$ By Lemma \ref{lemma:durett} the Markov chain $S_0, \cdots, S_n, \cdots
$ is transient. \ep \medbreak

In order to compare the dynamics $\mathcal{F}$ on the layers $L_0, L_1, \cdots,
L_n, \cdots$ with the Markov chain dynamics on the states $1, 2, \cdots, n, \cdots$
for any initial distribution $\rho$ on $\mathbb{T}^2$,
define $\mu=(\mathcal{F}_*\rho)|_{(\mathbb{T}^2 \setminus [-a,a]\times \{0\})}$ and $\nu$ on the Markov chain states as follows:
$\nu[i]=\mu(L_i)$, where $\nu[i]$ is the measure at the state $i$.
We are going to show that after the same number of iteration for both the
dynamics $\mathcal{F}$ on the layers and the Markov chain dynamics, the total
measure on the states of the Markov chain with $i \geq k$  is always
less or equal the total measure on the layers $\{ L_i, i \geq k \}$.
Denote by $(\mathcal{S}^n_*\nu)[i]$ the measure of the state $i$ at the $n^{th}$
iteration of the Markov chain.

\begin{lemma} \label{lemma:mc vs dynamics}
$$\forall n,\forall k \sum_{i=k}^{\infty}(\mathcal{S}_*^n\nu)[i] \leq \sum_{i=k}^{\infty}(\mathcal{F}^n_*\mu)(L_i)$$
\end{lemma}

\proc{Proof.} We are going to prove this lemma by induction on $n$ (for all $k$
simultaneously). \\
By the definition of $\nu$, the base of induction
$(n=0)$ is true: $\forall k \sum_{i=k}^{\infty}\nu[i] \leq
\sum_{i=k}^{\infty}\mu(L_i)$.

Suppose for $n=N$ the following holds:
$$ \forall k \ \ \sum_{i=k}^{\infty}(\mathcal{S}_*^N\nu)[i] \leq \sum_{i=k}^{\infty}(\mathcal{F}^N_*\mu)(L_i)$$
We would like to show that the same statement holds for $n=N+1$ given the statement for $n=N$, i.e.:
$$ \forall k \ \ \sum_{i=k}^{\infty}(\mathcal{S}_*^{N+1}\nu)[i] \leq \sum_{i=k}^{\infty}(\mathcal{F}^{N+1}_*\mu)(L_i)$$

The result is obvious for $k=0$: the total measure is always the same.

For $k>1$:\\
The pushed forward measure on the layers from $L_k$ through $L_\infty$
comes from several sources: all the measure from layers $L_{k+1}$ through
$L_\infty$, at least $\frac{2}{3}$ of the measure from the layer $L_k$, at
least $\frac{2}{3}$ of the measure from the layer $L_{k-1}$, and possibly
some extra from the layers $L_0, \cdots L_{k-2}$. We will only count the
first three sources with the $\geq$ sign.

$$\sum_{i=k}^{\infty}(\mathcal{F}^{N+1}_*\mu)(L_i) \geq \sum_{i=k+1}^{\infty}(\mathcal{F}^N_*\mu)(L_i)+\frac{2}{3}(\mathcal{F}^N_*\mu)(L_k)+\frac{2}{3}(\mathcal{F}^N_*\mu)(L_{k-1})$$

For the Markov chain we get exactly the first three sources with \textquotedblleft at
least $\frac{2}{3}$ of the measure" replaced by  \textquotedblleft exactly
$\frac{2}{3}$".

$$\sum_{i=k}^{\infty}(\mathcal{S}_*^{N+1}\nu)[i] = \sum_{i=k+1}^{\infty}(\mathcal{S}_*^N\nu)[i] + \frac{2}{3}(\mathcal{S}_*^N)\nu[k]+\frac{2}{3}(\mathcal{S}_*^N\nu)[k-1]$$

Applying the induction assumption for $k-1$ and $k+1$:

$$\sum_{i=k}^{\infty}(\mathcal{S}_*^{N+1}\nu)[i] = \sum_{i=k+1}^{\infty}(\mathcal{S}_*^N\nu)[i] + \frac{2}{3}(\mathcal{S}_*^N\nu)[k]+\frac{2}{3}(\mathcal{S}_*^N\nu)[k-1] = $$
$$=\frac{2}{3}\sum_{i=k-1}^{\infty}(\mathcal{S}_*^N\nu)[i] + \frac{1}{3}\sum_{k+1}^\infty (\mathcal{S}_*^N\nu)[i] \leq $$
$$
\frac{2}{3}\sum_{i=k-1}^{\infty}(\mathcal{F}^N_*\mu)(L_i)+\frac{1}{3}\sum_{i=k+1}^{\infty}(\mathcal{F}^N_*\mu)(L_i)
\leq $$
$$
\sum_{i=k+1}^{\infty}(\mathcal{F}^N_*\mu)(L_i)+\frac{2}{3}(\mathcal{F}^N_*\mu)(L_k)+\frac{2}{3}(\mathcal{F}^N_*\mu)(L_{k-1})
\leq \sum_{i=k}^{\infty}(\mathcal{F}^{N+1}_*\mu)(L_i) $$

For $k=1$: \\
The situation is just slightly different here: the pushed forward
measure on the layers from $L_1$ through $L_\infty$ comes from three
sources, two as before, while the amount of measure that comes
from the layer $L_0$ (or Markov chain state $0$) has a different coefficient.
Therefore:
$$\sum_{i=1}^{\infty}(\mathcal{F}^{N+1}_*\mu)(L_i) \geq \sum_{i=2}^{\infty}(\mathcal{F}^N_*\mu)(L_i)+\frac{2}{3}(\mathcal{F}^N_*\mu)(L_1)+q(\mathcal{F}^N_*\mu)(L_0)$$

$$\sum_{i=1}^{\infty}(\mathcal{S}_*^{N+1}\nu)[i] = \sum_{i=2}^{\infty}(\mathcal{S}_*^N\nu)[i] + \frac{2}{3}(\mathcal{S}_*^N\nu)[1]+q(\mathcal{S}_*^N \nu)[0]$$

Applying the induction assumption for $k= 0,1,2$:

$$\sum_{i=1}^{\infty}(\mathcal{S}_*^{N+1}\nu)[i] = \sum_{i=2}^{\infty}(\mathcal{S}_*^N\nu)[i] + \frac{2}{3}(\mathcal{S}_*^N\nu)[1]+q(\mathcal{S}_*^N\nu)[0]$$
$$ \frac{1}{3}\sum_{i=2}^{\infty}(\mathcal{S}_*^N\nu)[i] +
(\frac{2}{3}-q)\sum_{i=1}^{\infty}(\mathcal{S}_*^N\nu)[i] +
q\sum_{i=0}^{\infty}(\mathcal{S}_*^N\nu)[i] \leq $$
$$\frac{1}{3}\sum_{i=2}^{\infty}(\mathcal{F}^N_*\mu)(L_i)+
(\frac{2}{3}-q)\sum_{i=1}^{\infty}(\mathcal{F}^N_*\mu)(L_i) +
q\sum_{i=0}^{\infty}(\mathcal{F}^N_*\mu)(L_i) \leq $$
$$\sum_{i=2}^{\infty}(\mathcal{F}^N_*\mu)(L_i)+\frac{2}{3}(\mathcal{F}^N_*\mu)(L_1)+q(\mathcal{F}^N_*\mu)(L_0)
\leq \sum_{i=1}^{\infty}(\mathcal{F}^{N+1}_*\mu)(L_i) $$

\ep \medbreak

\proc{Proof of Prop. \ref{prop:local}}

Assume $\rho$ is $\mathcal{F}$-invariant probability measure on
$\mathbb{T}^2 \setminus ([-a,a]\times \{0\})$ and define the initial
measure for the Markov chain to be $\tau$ such that
$\tau[i]=\rho(L_i)$. As we push both of them forward, $\tau$ escapes
to $\infty$ by chain transience, implying that the measure of any
finite collection of states goes to 0. On the other hand, if $\rho$
stays invariant, it fails to satisfy the Lemma \ref{lemma:mc vs dynamics}. Therefore, we can conclude that there does not exist a $\mathcal{F}$-invariant
probability measure on $\mathbb{T}^2 \setminus ([-a,a]\times \{0\})$,
implying that all the invariant measures in the system must be
supported on $[-a,a] \times \{0\}$. The $\mathcal{F}$-invariant singular
measure on $[-a,a] \times \{0\}$ is unique by the 1-dimensional contraction example 1.1.

If we start with any Borel probability measure $\rho$ on $\mathbb{T}^2$,
$\mathcal{F}^n_*\rho$ must have a limit point by compactness. Define $\mu=(\mathcal{F}_*\rho)|_{(\mathbb{T}^2 \setminus [-a,a]\times \{0\})}$ and the initial
measure for the Markov chain $\tau$ such that $\tau[i]=\mu(L_i)$. Lemma \ref{lemma:mc vs dynamics} and Markov chain transience imply that all the limit points of $\mathcal{F}^n_*\mu$ and of $\mathcal{F}^n_*\rho$ are singular measures supported on $[-a,a] \times \{0\}$. \ep \medbreak

\section{Global Analysis}

Let $f:\mathbb{T}^2 \to \mathbb{T}^2$ be the Cat Map generated by the matrix  $\begin{pmatrix} 2 & 1 \\ 1 & 1 \\
\end{pmatrix}$.

Define $W \subset \mathbb{T}^2$ as before to be a neighborhood around the fixed
point bounded by the curves $x=\pm a$ and $y=\pm b\exp\{-\frac{x^2}{2\sigma^2}\}$ for
some $a,b, \sigma >0$. The local dynamics $f(x,y)=(\mu x, \lambda y)$ inside $W$ is exactly as we defined before with $\lambda$ and $\mu$ being the two
eigenvalues of the matrix.

Let $U$ and $U'$, $W \subset U' \subset U$, be two neighborhoods of $W$ bounded by the Gaussian-shaped curves $y=\pm(b+\beta)\exp\{-\frac{x^2}{2\sigma^2}\}$ and line segments $x=\pm (a+\beta)$ and $x=\pm (a+\frac{1}{2}\beta)$ respectively for some $\beta>0$. $C^\infty$ extend the vector field in $W$ such that it forms small constant angle $\varphi$ with the unstable direction outside the neighborhood $U'$. Denote the resulting vector field by $V$.

In order to produce the same effect as in the previous section, we need the parameters  to satisfy provisional assumptions similar to (A1) and (A2). The following assumptions will achieve the goal. Let
$\epsilon$, $a=\epsilon/(1-\mu)$, $b$, $\beta$, $\sigma$, and $\varphi$ be such
that:
\begin{itemize}
\item (B1) $\epsilon>3(\sigma
\sqrt{\frac{4ln\lambda}{1-\mu^2}} + b)$
\item (B2) $\beta<\frac{\epsilon}{7}$ and $\beta \leq (\lambda-1)b$
\item (B3) $\frac{1}{4}\beta>\frac{\epsilon \sin(\varphi)}{1-\mu}+(b+\beta)\exp\{-\frac{(a+\beta)^2}{2\sigma^2}\}\tan(\varphi)$.
\end{itemize}

\begin{theo} \label{thm:cat map}
Let $f$ be the Cat Map, $V$ and $\epsilon$ defined as above with parameters satisfying (B1),(B2) and (B3). Then
\begin{enumerate}
\item The only $\mathcal{F}$-invariant measure is the unique singular measure $\mu$ supported on the segment $[-a,a] \times \{0\} = [-\frac{\epsilon}{1-\mu}, \frac{\epsilon}{1-\mu}] \times \{0\}$ of the local stable manifold of the fixed point in the $E^s-E^u$-coordinate frame.
\item If $\nu$ is any Borel probability measure on $\mathbb{T}^2$, then $\mathcal{F}^n_* \nu$ has all its limit points supported on $[-a,a] \times \{0\}$. In particular, $\frac{1}{n}\sum_{i=0}^{n-1} \mathcal{F}_*^i \nu \to \mu$.
\end{enumerate}
\end{theo}

First, we are going to show that for the Cat Map the return rates from
$(\mathbb{T}^2 \setminus W)$ to W are similar to the assumption (A2)
from the previous section. Then we are going to model the dynamics with an appropriate transient countable state Markov chain as in section 2.

\subsection{Return Rates For The Cat Map}

\begin{prop} \label{prop:return rates}
Let $f$ be the Cat Map, $V$ and $\epsilon$ defined as above with parameters satisfying (B1),(B2) and (B3). Then there exists $N$ and $\eta>0$ such that for any finite measure $\nu$ on $\mathbb{T}^2$, at least $\eta$ fraction of $\nu|_{W^c}$ gets to $W$ in $N$ steps when pushed forward by the perturbed dynamics $\mathcal{F}$. I.e. if $\mu$ is such that $\mu(A)=\nu(A \cap W^c)$, then $\mathcal{F}^N_*\mu(W) \geq \eta \mu(W^c)$.
\end{prop}

Let $\nu$ be any finite measure supported in $W^c$. If we
push it forward by the function $f$ once, we can divide it into three
parts:
\begin{enumerate}
\item $(f_*\nu)|_W$,
\item $(f_*\nu)|_{U \setminus W}$, and
\item $(f_*\nu)|_{U^c}$
\end{enumerate}

We are going to follow what happens when we perturb each
of these three parts separately and then push the corresponding
measures forward by the perturbed dynamics $\mathcal{F}$ some more
times. If we can pick a single $N$ for all of these measures such that after
the total of $N$ push forwards by the perturbed dynamics, certain
fraction $\eta$ of the initial measure stays in $W$, that will prove Proposition \ref{prop:return rates}.

\begin{lemma} \label{lemma:W returns}
A lower bound estimate for the fraction of the measure that ends up in
W after the perturbation of $(f_*\nu)|_W$ followed by $(n-1)$ push
forwards by $\mathcal{F}$, $n \geq 1$, is $(\frac{1}{2})(\frac{2}{3})^{n-1}$.
\end{lemma}

\proc{Proof.}
After the initial perturbation at most $\frac{1}{2}$-fraction of the
measure can \textquotedblleft escape" through the ends of $W$ since the vector field in $W$ is Gaussian-shaped and $\epsilon < 2a$. For all the future
iterations, the part of the measure that remains inside $W$ after being pushed forward by $f$ stays in $W$ after the perturbation as well. Since
the dynamics inside $W$ is exactly as described in section 2 and (A1)=(B1), we can apply Lemma \ref{lemma:2/3}. Thus at most
$\frac{1}{3}$ of the measure can \textquotedblleft escape" to $W^c$ in any single step and
at least $\frac{2}{3}$-fraction of the measure supported in $W$
stays in $W$ after being pushed forward by $\mathcal{F}$. Therefore, $(\frac{1}{2})(\frac{2}{3})^{n-1}$-fraction of $(f_*\nu)|_W$ stays in $W$ after the perturbation followed by $n-1$ push forwards by $\mathcal{F}$. \ep \medbreak

\begin{lemma} \label{lemma:U-W returns}
A lower bound estimate for the fraction of the measure that ends up in
$W$ after the perturbation of $(f_*\nu)|_{U \setminus W}$ followed by
$(n-1)$ push forwards by $\mathcal{F}$, $n \geq 2$, is
$(\frac{d}{2\epsilon \lambda})(\frac{2}{3})^{n-2}$, where $d=b \exp\{-\frac{a^2}{2\sigma^2}\}$.
\end{lemma}

\proc{Proof}

\begin{figure}[center]
  \centering{\scalebox{0.7}{\includegraphics{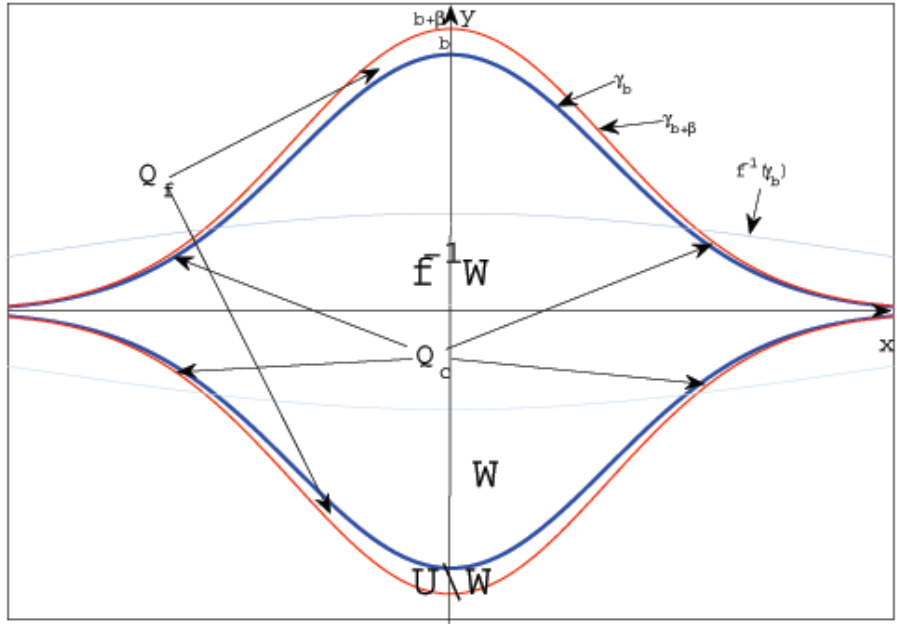}}}\\
  \caption{Neighborhood U with $Q_f$ and $Q_c$ parts}\label{fig:U}
\end{figure}

We are going to split the neighborhood $U \setminus W$ into two parts and argue about each separately. Let $\pm \gamma_b$ and $\pm \gamma_{b+\beta}$ be the Gaussian-shaped curves given by equations
$y=\pm be^{-\frac{x^2}{2\sigma^2}}$ and $y=\pm (b+\beta)\exp\{-\frac{x^2}{2\sigma^2}\}$ respectively. The curves $f^{-1}(\pm \gamma_b)$ split $U \setminus W$ into two parts: $Q_c = (U \setminus W) \cap f^{-1}W$ and $Q_f = (U \setminus W) \setminus f^{-1}W$, where \textquotedblleft$c$" stands for \textquotedblleft closer to the $x$-axis" and \textquotedblleft$f$" for \textquotedblleft farther from the $x$-axis". See Figure \ref{fig:U}.

\textbf{$Q_f$ part:}
One can guess that the worst case for the $Q_f$ part occurs when $(f_*\nu)|_{U \setminus W}$ is a point measure at $(x,(b+\beta)\exp\{-\frac{x^2}{2\sigma^2}\})$ for $x=-\sigma \sqrt{\frac{4ln(\lambda)}{1-\mu^2}}$. In this situation half of the measure \textquotedblleft escapes" from $U$ after the perturbation (by the nature of the vector field) and a big portion of the measure that stays \textquotedblleft close" to the $y$-axis after the perturbation will \textquotedblleft escape" $U$ on the next push forward by $f$. The goal is to estimate the amount of measure that gets to $f^{-1}W$ after the first perturbation and therefore does not leave $W$ at the next iteration.

We assumed in (B1) that $\epsilon>3(\sigma \sqrt{\frac{4ln\lambda}{1-\mu^2}} +
b)$ and in (B2) that $\beta < \frac{\epsilon}{7}$ and $\frac{b+\beta}{b}<\lambda$. The following statements hold:
\begin{itemize}

\item The curves
$\gamma_{b+\beta}$ and
$f^{-1}\gamma_b$ intersect with $x$-coordinate $\pm
\sigma \sqrt{\frac{2ln(\lambda(b+\beta)/b)}{1-\mu^2}}$ and by (B2) $\sigma \sqrt{\frac{2ln(\lambda(b+\beta)/b)}{1-\mu^2}} <$ $\sigma \sqrt{\frac{4ln\lambda}{1-\mu^2}}$

\item The distance along $\gamma_{b+\beta}$ from $(0,(b+\beta))$ to $f^{-1}W$ is at most $\sigma \sqrt{\frac{2ln(\lambda(b+\beta)/b)}{1-\mu^2}} +
b+\beta$ by an argument similar to the proof of the Lemma \ref{lemma:epsilon}; and by (B1) and (B2), $\sigma \sqrt{\frac{2ln(\lambda(b+\beta)/b)}{1-\mu^2}} +
b+\beta <$ $\sigma \sqrt{\frac{4ln\lambda}{1-\mu^2}} +b + \beta <$ $\frac{\epsilon}{3}+\frac{\epsilon}{7} = \frac{10\epsilon}{21}$

\item The support of the push forward of $\nu|_{Q_f}$ under the perturbation is located within the strip bounded by $x=\pm (\sigma \sqrt{\frac{2ln(\lambda(b+\beta)/b)}{1-\mu^2}}+\epsilon)$ and by (B1) and (B2), $\sigma \sqrt{\frac{2ln(\lambda(b+\beta)/b)}{1-\mu^2}}+\epsilon <$ $\sigma \sqrt{\frac{4ln\lambda}{1-\mu^2}}+\epsilon <$ $\frac{1}{3}\epsilon +
\epsilon=\frac{4}{3}\epsilon<\frac{\epsilon}{1-\mu}=a$.
\end{itemize}

Therefore, when $(f_*\nu)|_{Q_f}$ is perturbed, at least $\frac{2\epsilon - \epsilon - 2\frac{10}{21}\epsilon}{2\epsilon}=\frac{1}{42}$-fraction of $(f_*\nu)|_{Q_f}$ gets into $f^{-1}W$; same fraction clearly stays in $W$ following a push forward by $f$. When perturbed again, the measure that ended up in $f^{-1}W$ after the first perturbation does not leave $W$ since
$$\mu(\sigma \sqrt{\frac{2ln(\lambda(b+\beta)/b)}{1-\mu^2}}+\epsilon) +
\epsilon < \mu\sigma\sqrt{\frac{4ln\lambda}{1-\mu^2}}+
\mu\epsilon + \epsilon < \frac{1}{3}\epsilon\mu + \mu\epsilon +
\epsilon=$$
$$(\frac{4}{3}\mu + 1)\epsilon<\frac{\epsilon}{1-\mu}=a.$$

After the next $(n-2)$ push forwards under $\mathcal{F}$, at least
$(\frac{2}{3})^{n-2}$-fraction of the measure will stay in $W$, so
together it ensures that at least $(\frac{1}{42})(\frac{2}{3})^{n-2}$-fraction of $(f_*\nu)|_{Q_f}$ will end up in $W$ after $n$ steps.

\textbf{$Q_c$ part:}

\begin{figure}[center]
  \centering{\scalebox{0.7}{\includegraphics{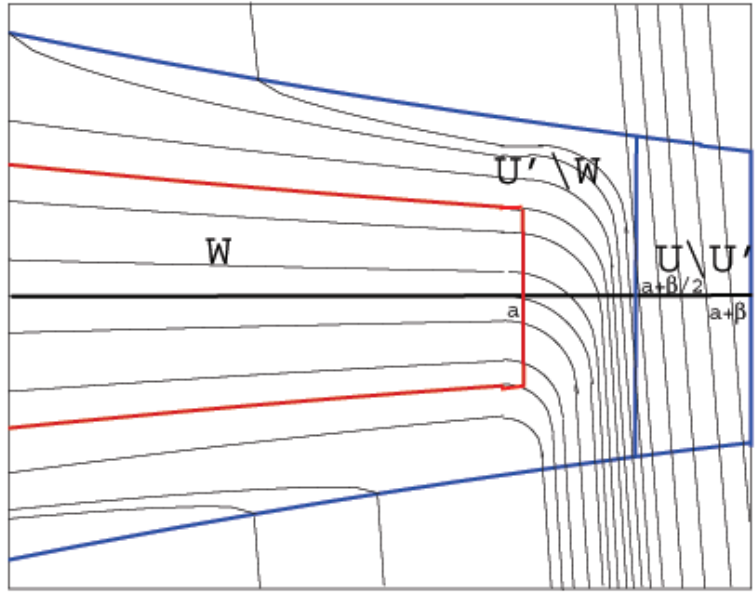}}}\\
  \caption{Vector field at the ends}\label{Figure}
\end{figure}

Let $d'=(b+\beta)\exp\{-\frac{(a+\beta)^2}{2\sigma^2}\}$.

In this case the worst estimate occurs for the point measure located at one of the \textquotedblleft corners" of $U$, e.g. at $(a+\beta,d')$ if $\varphi < 0$ as in Fig \ref{Figure}.

After the first perturbation, we can guarantee that at least
$\frac{2d'}{2\epsilon \cos(\varphi)}$-fraction of $(f_*\nu)|_{Q_c}$ stays in $U \cup (([a+\beta,a+\beta+2d'\tan(\varphi)] \cup [-a-\beta,-a-\beta-2d'\tan(\varphi)])\times[-d',d'])$. When
we push that part of the measure forward by $f$, at least $\frac{2d \cos(\varphi)}{2\lambda d' \cos(\varphi)}=\frac{d}{\lambda d'}$-fraction of it will
end up in $W$. After another perturbation, the fraction of the measure that remains in $W$ is at least
$$\epsilon+(a-\mu(a+\beta+d'\tan(\varphi)))=\epsilon+a(1-\mu)-\mu\beta-d'\tan(\varphi)>$$
$$>2\epsilon-\mu\beta-\mu \frac{1}{4}\beta >2\epsilon - \frac{5}{4}\frac{1}{7}\mu \epsilon>\epsilon(2-\frac{5}{28}\mu)>\epsilon,$$

by assumptions (B2) and (B3). Therefore, we can ensure that at least $\frac{2d'}{2\epsilon \cos(\varphi)} \frac{d}{\lambda d'} \frac{\epsilon}{2\epsilon}=$ $\frac{d}{2\epsilon \lambda \cos(\varphi)}>\frac{d}{2 \epsilon \lambda}$-fraction of $(f_*\nu)|_{Q_c}$ ends up in $W$ after the perturbation and one push forward under $\mathcal{F}$; and at least $(\frac{d}{2\epsilon \lambda })(\frac{2}{3})^{n-2}$-fraction of $(f_*\nu)|_{Q_c}$ will end up in $W$ after $n$ steps, $n \geq 2$.

Assumption (B1) states that $\epsilon>3(\sigma \sqrt{\frac{4ln\lambda}{1-\mu^2}} + b)$, thus $\sigma^2 < \frac{ \epsilon^2 (1-\mu^2)}{9 \cdot 4 ln\lambda}$ and

$$\exp\{-\frac{a^2}{2\sigma^2}\}=\exp\{-\frac{\epsilon^2}{2(1-\mu)^2\sigma^2}\}<
\exp\{-\frac{36 \epsilon^2 ln \lambda}{2(1-\mu)^2 \epsilon^2 (1-\mu^2)}\}=$$
$$=\exp\{-\frac{18 ln \lambda}{2(1-\mu)^2(1-\mu^2)}\}<10^{-23}<\frac{1}{100}$$

Therefore
$$\frac{d}{2\epsilon \lambda}=\frac{be^{-\frac{a^2}{2\sigma^2}}}{2\epsilon \lambda }<\frac{b}{200 \epsilon \lambda}<\frac{\epsilon}{3 \cdot 200 \epsilon \lambda}<\frac{1}{42}.$$

So a lower bound estimate for the fraction of the measure that ends up in
$W$ after the perturbation of $(f_*\nu)|_{U \setminus W}$ followed by
$(n-1)$ push forwards by $\mathcal{F}$, $n \geq 2$, is
$(\frac{d}{2\epsilon \lambda})(\frac{2}{3})^{n-2}$. \ep \medbreak

\begin{lemma} \label{lemma:U^c returns}
There exist $N$ and $\eta>0$ (both do not depend on $\nu$) such that a lower bound estimate for the fraction of the measure that ends up in $W$ after the perturbation of $(f_*\nu)|_{U^c}$ followed by $(N-1)$ push forwards by $\mathcal{F}$ is $\eta$.
\end{lemma}

Instead of dealing with the push forwards of the measure $\nu$, we are going to simplify the discussion by dealing with the push forwards of the measures $\delta_x$ with $f(x) \in U^c$. To prove Lemma \ref{lemma:U^c returns} it is enough to show that there exist $N$ and $\eta>0$ such that for any $x \in \mathbb{T}^2$ with $f(x) \in U^c$, $(\mathcal{F}^N_* \delta_x)(W) \geq \eta$.

To simplify the discussion even further, we are going to start with proving Lemma \ref{lemma:U^c returns} for the random perturbation of the Cat Map $f$ along the vector field $V'$ that agrees with $V$ in $W$, but is everywhere parallel to the unstable direction in $(U')^c$ (with $C^\infty$ extension in $U'\setminus W$). Let $\mathcal{F}'$ denote the perturbation of the Cat Map $f$ along $V'$ with the perturbation size $\epsilon$ and $\mathcal{P}'$ denote the perturbation of the identity map along $V'$ with the perturbation size $\epsilon$.

\begin{lemma} \label{lemma: V' return rates}
There exist $N' \geq 3$ such that for any $x \in \mathbb{T}^2$ with $f(x) \in U^c$ $$((\mathcal{F}')_*^{N'} \delta_x)(W) \geq \kappa_{N'-2} (\frac{d}{2\epsilon \lambda}) (\frac{2}{3})^{N'-2},$$
where $\kappa_n$ does not depend on $x$.
\end{lemma}

\proc{Proof.}
$f_*(\mathcal{F}_*'\delta_x)$ is supported on a $2\lambda \epsilon$-long interval parallel to the unstable direction. If it fully crosses $U_0=[-a-\beta,a+\beta] \times [-d',d']$, then at least $\kappa_0=\frac{d'}{\lambda \epsilon}$ fraction of the $f_*(\mathcal{F}_*'\delta_x)$-measure gets to $U$. By fully crossing we mean that the interval fully goes through the rectangle. i.e., in this case, their intersection is of length $2d'$.

Whether or not supp$[f_*(\mathcal{F}_*'\delta_x)]$ fully crosses $U_0$, if $[f_*(\mathcal{F}_*'\delta_x)](U) \geq \kappa_0$, then by the estimates from Lemmas \ref{lemma:W returns} and \ref{lemma:U-W returns}
$$((\mathcal{F}')^{N'}_* \delta_x)(W) \geq \kappa_{0} (\frac{d}{2\epsilon \lambda}) (\frac{2}{3})^{N'-2}$$ for any $N' \geq 3$. We will pick appropriate $N'$ that works for all $x$ with $f(x) \in U^c$ later in this proof.

If $[f_*(\mathcal{F}_*'\delta_x)]< \kappa_0$, supp$[f_*(\mathcal{F}_*'\delta_x)]$ does not fully cross $U_0$, implying that supp$[f_*(\mathcal{F}_*'\delta_x)] \cap U^c$ consists of one piece. Let $I_1'=$supp$[f_*(\mathcal{F}_*'\delta_x)] \cap U^c$.
Then $I_1'$ is an interval parallel to the unstable direction with $\lambda \epsilon \leq |I_1'| \leq 2 \lambda \epsilon$.

Either $f(I_1')$ fully crosses $U_0$ or $f(I_1') \cap U^c$ consists of a single piece.
While $f(I_n') \cap U^c$ consists of a single piece, define inductively $I_{n+1}'=f(I_n') \cap U^c$. Clearly $I_n' \subset$ supp$[f_*((\mathcal{F}')^n_*\delta_x)]$. We would like to estimate the growth of the length $|I_n'|$ before $f(I_n')$ fully crosses $U_0$. If $f(I_n') \cap (U)^c$ consists of a single piece, the intersection of $f(I_n')$ with $U^c$ \textquotedblleft chops off" at most $2(b+\beta+d')$ length from $f(I_n')$.

Using the estimate from the proof of Lemma \ref{lemma:U-W returns}, (B1), and (B2), $$d'=(b+\beta)\exp\{-\frac{(a+\beta)^2}{2\sigma^2}\}<(b+\beta)\exp\{-\frac{a^2}{2\sigma^2}\} < (b+\beta)\frac{1}{100}< \frac{1}{100} (\frac{\epsilon}{3}+\frac{\epsilon}{7})=\frac{\epsilon}{210}$$

and
$$2(b + \beta + d') < 2(\frac{\epsilon}{3}+\frac{\epsilon}{7}+ \frac{\epsilon}{210}) < 2\frac{\epsilon}{2}=\epsilon.$$

We conclude that while $f(I_n') \cap U^c$ consists of a single piece, the length of $I_n'$ grows exponentially with the rate at least $(\lambda-1) > 1.6$. Indeed, $|I_{n+1}'| \geq \lambda |I_n'| - \epsilon \geq \lambda|I_n'|-|I_1'| \geq (\lambda-1)|I_n'|$.

The following lemma states that if we take an interval parallel to the unstable direction long enough, it wraps around the torus densely enough to intersect any piece of stable manifold of a predetermined length.

\begin{lemma} \label{lemma:wrapping around the torus}
For any $\alpha>0$, there exists $l>0$ such that, if $I$ is a line segment of length at least $l$ parallel to the unstable direction of the Cat Map, then $I$ crosses any interval parallel to the stable direction of length at least $\alpha$. Moreover, if $I$ forms a constant angle $\varphi$, $-\pi/2 < \varphi < \pi/2$, with the unstable direction, the same result holds.
\end{lemma}

\proc{Proof.}
In the regular coordinates, the eigenvectors for the Cat Map corresponding to $\lambda$ and $\mu$ can be taken $(1,\tau)$ and $(-\tau,1)$ respectively, where $\tau = \frac{\sqrt{5}-1}{2}$. Assume we take the unstable manifold of the fixed point $(0,0)$ and look at the coordinates at which it intersects the $y$-axis. If we move in the positive sense with respect to the $y$-axis and the eigenvector $(1,\tau)$, the intersection points will be $0$, $\tau$, $2\tau \mod 1$, $\cdots$, which represent precisely the rotations of a circle with the rotation number $\tau$. By the properties of irrational rotations of a circle, there exists $n$ such that any interval on $y$-axis of length at least $\alpha \sqrt{1+\mu}$ contains at least one of any $n$ subsequent intersection points of the unstable manifold with the $y$-axis, i.e. it crosses any piece of unstable manifold of length at least $l=\frac{n}{\sqrt{1+\mu}}$. Therefore any given interval of length $\alpha$ parallel to the stable manifold contains at least one intersection point with an unstable manifold of length $\frac{n}{\sqrt{1+\mu}}$, which proves Lemma \ref{lemma:wrapping around the torus} for the situation when $I$ is parallel to the unstable direction. If $I$ forms a constant angle $\varphi$, $-\pi/2 < \varphi < \pi/2$, with the unstable direction, the distances between the intersection points on a piece of stable manifold are unchanged, implying that Lemma \ref{lemma:wrapping around the torus} holds when an interval forms an angle $\varphi$ with the unstable direction with a choice of $l=\frac{n}{\sqrt{1+\mu} \cos(\varphi)}$. \ep \medbreak

For the rest of this paper, we are going to work with the stable-unstable coordinate frame with $(0,0)$ being the fixed point of the Cat Map.

Let $N'$ be such that given $\alpha = 2a + 2\beta$, $l=\lambda^2 \epsilon (\lambda-1)^{N'-2}-2d'$ works for Lemma \ref{lemma:wrapping around the torus}. Then there exists $1 \leq k \leq N'-2$ such that $f(I_k')$ fully crosses $U_0$. We would like to establish that there exists a lower bound on the measure supported on $f(I_k') \cap U_0$ that does not depend on $x$.

Suppose $f(I_k')$ fully crosses $U_0$, $k \geq 1$. Let $I_1'=f^{-k}(f(I_k) \cap U_0)$ be a subinterval of $I_1$ that maps to $(f(I_k) \cap U_0)$ under $k$ iterations of the Cat Map. The length of $I_1$ is $|I_1'|=\frac{2d'}{2\epsilon \lambda^{k}}$ and it supports $\frac{|I_1'|}{2\epsilon}=\frac{2d'}{2 \lambda \epsilon \lambda^{k}}$-fraction of $f_*(\mathcal{F}_*' \delta_x)$-measure. When the perturbation $\mathcal{P}'$ occurs, the measure supported on $I_1'$ gets \textquotedblleft smeared"; and the worst estimate for the measure that remains in $I_1'$ after the perturbations occurs when we do not count the measure that gets to $I_1'$ from nearby parts of $I_1$. Therefore at least $\frac{|I_1'|}{2\epsilon}$-fraction of the measure stays in $I_1'$ after the perturbation. Define inductively $I_j'=f(I_{j-1}')$, $2 \leq j \leq k$. Then the $f$ push forward of any measure supported on $I_{j-1}'$ is supported on $I_j'$ and when a measure supported on $I_j'$ is perturbed, at least $\frac{|I_j'|}{2\epsilon}$-fraction of it stays in $I_j'$. Therefore, the amount of measure supported on $I_k'=(f(I_k) \cap U_0)$ is at least
$$\kappa_k=\frac{d'}{\epsilon \lambda^{k+1}} \frac{d'}{\epsilon \lambda^{k}} \cdots \frac{d'}{\epsilon \lambda^2}=\frac{(d')^k}{\epsilon^k \lambda^{\frac{k(k+3)}{2}}}.$$

Clearly $\kappa_n > \kappa_{n'}$ for $n<n'$. We conclude that, there exists $k \leq N'-2$ such that $[f_*((\mathcal{F}')^k_* \delta_x)](U) \geq \kappa_{N'-2}$. Note that $\kappa_{N'-2}$-fraction of measure may end up in $U$ before the support of the pushed forward measure fully crosses $U_0$. By the estimates from Lemmas \ref{lemma:W returns} and \ref{lemma:U-W returns}, $$((\mathcal{F}')^{N'}_* \delta_x)(W) \geq \kappa_{N'-2} (\frac{d}{2\epsilon \lambda}) (\frac{2}{3})^{N'-2}.$$

\ep \medbreak

This completes the proof of Lemma \ref{lemma:U^c returns} for the random perturbation of the Cat Map that occurs along the vector field $V'$ parallel to the unstable direction in $(U')^c$. Now we are ready to deal with vector field $V$ forming an angle $\varphi$ with the unstable direction in $(U')^c$.

\begin{lemma} \label{lemma: V return rates}
There exist $N$ such that for any $x \in \mathbb{T}^2$ with $f(x) \in U^c$, $$((\mathcal{F})_*^{N} \delta_x)(W) \geq \kappa_{(\varphi, N-2)} (\frac{d}{2\epsilon \lambda}) (\frac{2}{3})^{N-2},$$ where $\kappa_{(\varphi,n)}$ does not depend on $x$.
\end{lemma}

\proc{Proof.}
$f_*(\mathcal{F}_*\delta_x)$ is supported on an interval $I$ that forms angle $\varphi'$ with the unstable direction with $\tan(\varphi')=\frac{\mu}{\lambda}\tan(\varphi)$ and the length of its projection to the unstable direction is $|I|_u=|I|\cos(\varphi')=2\epsilon \cos(\varphi) \lambda$. If $I$ fully crosses $U_0=[-a-\beta,a+\beta] \times [-d',d']$, then at least $\kappa_{(\varphi,0)}=\frac{d'}{\lambda \epsilon \cos(\varphi)}$ fraction of the $f_*(\mathcal{F}_*\delta_x)$-measure gets to $U$. By fully crossing we mean that the interval fully goes through the rectangle, i.e. in this case their intersection is of length $\frac{2d'}{\cos(\varphi')}$. Then by the estimates from Lemmas \ref{lemma:W returns} and \ref{lemma:U-W returns},
$$(\mathcal{F}_*^{N} \delta_x)(W) \geq \kappa_{(\varphi, 0)} (\frac{d}{2\epsilon \lambda}) (\frac{2}{3})^{N-2},$$
for any $N \geq 3$.

If $[f_*(\mathcal{F}_* \delta_x)](U) < \kappa_{(\varphi,0)}$, define $I_1 =$ supp$[f_*(\mathcal{F}_* \delta_x)] \cap U^c$.

From (B2) and (B3) we have
$$\frac{\epsilon}{28} > \frac{1}{4}\beta> \frac{\epsilon \sin(\varphi)}{1-\mu}+be^{-\frac{a^2}{2\sigma^2}}\tan(\varphi) > \frac{\epsilon \sin(\varphi)}{1-\mu}.$$
Thus \small
$$\sin(\varphi)<\frac{1-\mu}{28} \Leftrightarrow 1-\cos^2(\varphi)=\sin^2(\varphi)<\frac{(1-\mu)^2}{28^2} \Leftrightarrow \frac{1}{\cos(\varphi)}<\sqrt{\frac{28^2}{28^2-(1-\mu)^2}}<1.1$$
\normalsize and
$$ \frac{b+\beta+d'}{\cos(\varphi)} < 1.1\frac{\epsilon}{2}.$$

Given $I_n$, define $I_{n+1}$ to be the middle $(\lambda-1.1)|I_n|$-part of $f(I_n)$, i.e. we \textquotedblleft chop off" $0.55\lambda|I_n|$ from both sides of $f(I_n)$. By the above estimates we \textquotedblleft chop off" enough to ensure that, unless $f(I_n)$ fully crosses $U_0=[-a-\beta, a+\beta]\times[-d',d']$, $I_{n+1}$ does not intersect with $U$. Then $|I_n|$ grows with rate at least $(\lambda-1.1)>1.5$. Let $N$ be such that, given $\alpha=2a+2\frac{3}{4}\beta$, $l=\lambda^2 \epsilon (\lambda-1.1)^{N-2}-2d'/\cos(\varphi)$ works for Lemma \ref{lemma:wrapping around the torus} with interval forming an angle $\varphi$ with the unstable direction; clearly same $l$ works for any interval forming an angle $\leq \varphi$ with the unstable direction. Then for some $k \leq N-2$, $f(I_k)$ fully crosses $U_0$ intersecting simultaneously with $[-a-\frac{3}{4}\beta,a+\frac{3}{4}\beta]$.

In the following, we would like to talk about parallelograms of the following kind: Given an interval $J$ forming some angle with the unstable direction, let $R(J,w)$ be the parallelogram with sides parallel to $J$ of length $|J|$ and sides parallel to the stable direction of length $w$ such that $J$ passes through the middle of $R(J,w)$ and divides it into two halves.

We defined $I_n$ for all $n$. Let $R_n=R(I_n,2\frac{\epsilon \sin(\varphi)}{1-\mu})$. In assumption (B3) we chose $\frac{1}{4}\beta > \frac{\epsilon \sin(\varphi)}{1-\mu}+d'\tan(\varphi)$. This guarantees that if $f(I_n)$ intersects $[-a-\frac{3}{4}\beta,a+\frac{3}{4}\beta]$ and fully crosses $U_0=[-a-\beta,a+\beta]\times[-d',d']$, then $f(R_n)$ also fully crosses $U_0$.

Since we defined $I_{n}$ to be the middle $(\lambda-1.1)|I_{n-1}|$-part of $f(I_{n-1})$, by the above estimates, $R_n \subset U^c$. Let $\pi$ be any measure supported on $R_n$. Then supp$f_* \pi \subset R(f(I_n),2\frac{\mu \epsilon \sin(\varphi)}{1-\mu}))$. If we then perturb $f_* \pi$, the maximal \textquotedblleft width" on the support of the resulting measure is $2\frac{\epsilon \sin(\varphi)}{1-\mu}$, i.e. when \textquotedblleft smearing" occurs, the measure can only \textquotedblleft leak" through the sides of $R_{n+1}$ parallel to the stable direction.

It follows that if $f(I_n)$ does \emph{not} intersect $[-a-\frac{3}{4}\beta,a+\frac{3}{4}\beta]$, $R(f(I_n),2\frac{\epsilon \sin(\varphi)}{1-\mu}) \cap U' = \emptyset$. Therefore the perturbation occurs at constant angle $\varphi$ throughout $R(f(I_n),2\frac{\epsilon \sin(\varphi)}{1-\mu})$. Although we cannot guarantee that enough measure gets to $U$ at this step, we can ignore such a crossing and continue our iteration.

Suppose $f(I_k)$ intersects $[-a-\frac{3}{4}\beta,a+\frac{3}{4}\beta]$ and fully crosses $U_0=[-a-\beta,a+\beta]\times[-d',d']$ for some $1 \leq k \leq N-2$. Let $I_j'=f^{-k+j}(f(I_k) \cap U_0)$, $1 \leq j \leq k$ and let $R_j'=R(I_j', 2\frac{\mu \epsilon \sin(\varphi)}{1-\mu})$. Note that we chose the widths of $R_j$'s and $R_j'$'s such that
\begin{itemize}
\item if we perturb any measure supported on $R_j'$, the maximal \textquotedblleft width" of the resulting measure support is going to be the width of $R_n$, i.e. $ 2\frac{\epsilon \sin(\varphi)}{1-\mu}$, except possibly at the ends; and
\item if we push any measure supported on $R(I_j',2\frac{\epsilon \sin(\varphi)}{1-\mu})$ forward under $f$, the resulting measure is supported in $R_{j+1}'$.
\end{itemize}

It follows that the fraction of any measure supported on $R_j'$ that stays in $R(I_j',2\frac{\epsilon \sin(\varphi)}{1-\mu})$ after the perturbation $\mathcal{P}$ is $\frac{|R_j'|_u}{2\epsilon \cos(\varphi)}=\frac{d'}{\epsilon \cos(\varphi) \lambda^{k-j+1}}$, where $|R_j'|_u$ denotes the length the projection of $R_j'$ to the unstable direction. From here we compute a lower bound $\kappa_{(\varphi,k)}$ on the amount of $f_*(\mathcal{F}^k_*\delta_x)$-measure contained in $f(R_k) \cap U_0$:
$$\kappa_{(\varphi,k)}=\frac{d'}{\epsilon \cos(\varphi) \lambda^{k+1}} \frac{d'}{\epsilon \cos(\varphi) \lambda^{k}} \cdots \frac{d'}{\epsilon \cos(\varphi) \lambda^2}=\frac{(d')^k}{\epsilon^k \cos^k(\varphi) \lambda^{\frac{k(k+3)}{2}}}$$

With $N$ and $\kappa_{(\varphi,N-2)}$ defined as above, we are guaranteed that for some $k \leq N-2$, $f_*(\mathcal{F}^k_* \delta_x)(U) \geq \kappa_{(\varphi,N-2)} $. Using the estimates from Lemmas \ref{lemma:W returns} and \ref{lemma:U-W returns} we conclude that
$$(\mathcal{F}^N_* \delta_x)(W) \geq \kappa_{(\varphi, N-2)} (\frac{d}{2\epsilon \lambda}) (\frac{2}{3})^{N-2}.$$

This proves Lemma \ref{lemma: V return rates}, which implies Lemma \ref{lemma:U^c returns}.

Lemmas \ref{lemma:W returns}, \ref{lemma:U-W returns}, and \ref{lemma:U^c returns} imply Proposition \ref{prop:return rates}.
\ep \medbreak

\subsection{Markov chain modifications}

Now we need to adjust the our argument to fit the return rates that
happen within $N$ steps instead of happening every step as we assumed
before. For that we want to look only what happens at each $N^{th}$
iteration of our dynamics, $N$ is as in Proposition \ref{prop:return rates}. Consider the countable state Markov chain described as follows:

\begin{itemize}
\item For $\eta$ as in Proposition \ref{prop:return rates}
$$P(0,1)=\eta\ \ and\ \ P(0,0)=1-\eta; $$
\item for $k \geq N$:
$$P(k,k+N)=C^0_N (\frac{2}{3})^N;\
P(k,k+N-2)=C^1_N(\frac{2}{3})^{N-1}(\frac{1}{3}); \cdots ;$$
$$P(k,k+N-2s)=C^{s}_N(\frac{2}{3})^{N-s}(\frac{1}{3})^s; \cdots;
P(k,k-N)=C^N_N(\frac{1}{3})^N;$$
\item $k < N$ and $s$ be such that $1 \leq k+N-2s \leq 2$
$$P(k,k+N)=C^0_N(\frac{2}{3})^N;\
P(k,k+N-2)=C^1_N(\frac{2}{3}^{N-1}(\frac{1}{3}); \cdots ;$$
$$P(k,k+N-2s)=C^{s}_N(\frac{2}{3})^{N-s}(\frac{1}{3})^s,
P(k,0)=1-\sum_{i=0}^s P(k,k+N-2i).$$

\end{itemize}

Here $C^k_N$ stands for binomial coefficient \textquotedblleft $n$ choose $k$."

Note that the transition probabilities are exactly the same as if
we looked at the $N^{th}$ iteration of the Markov chain from section 2
except at the first $N$ states. By another application of Lemma \ref{lemma:durett} with $\varphi(x)=\frac{1}{2^x} \Rightarrow E_k\varphi(X_1) = \varphi(k)$ we conclude the this Markov chain is transient.

We would like to show now that the $N^{th}$ snapshots of the perturbed dynamics can be compared to the dynamics of the Markov chain. Let $\rho$ be any measure on $\mathbb{T}^2$ and $\mu=(\mathcal{F}_*\rho)|_{(\mathbb{T}^2 \setminus ([-a,a] \times \{0\}))}$. Define $\nu$ on the Markov chain states to be $\nu[i]=\mu(L_i)$. Then after $nN$ number of
iteration of the dynamics $\mathcal{F}$ and $n$ iterations of the Markov chain, the total measure on the Markov chain states $\{[i], i \geq
k \}$ is always less or equal the total measure on the layers $\{
L_i, i \geq k \}$.

\begin{lemma} \label{lemma:mc vs dynamics modified}
$$\forall n,\forall k \sum_{i=k}^{\infty}\mathcal{S}_*^n\nu[i] \leq \sum_{i=k}^{\infty}\mathcal{F}^{nN}_*\mu(L_i)$$
\end{lemma}

The proof of this Lemma is exactly the same as for Lemma \ref{lemma:mc vs dynamics}, by induction on $n$ for all $k$ simultaneously, except it involves more terms and thus omitted.

\proc{Proof of Theorem \ref{thm:cat map}.}
Suppose $\rho$ is $\mathcal{F}$-invariant measure on $\mathbb{T}^2
\setminus ([-a,a]\times \{0\})$  and define the initial measure for the
Markov chain to be $\tau$ such that $\tau[i]=\rho(L_i)$. As we
push both of them forward, $\tau$ escapes to $\infty$ by chain
transience, implying that the measure of any finite collection of
states goes to 0. On the other hand, if $\rho$ stays invariant, it
fails to satisfy the Lemma \ref{lemma:mc vs dynamics modified}. Therefore, we can conclude that
there does not exist a $\mathcal{F}$-invariant measure on $\mathbb{T}^2
\setminus ([-a,a]\times \{0\})$, the unique invariant measure in the system is supported on the segment of the local stable manifold of the fixed point $[-a,a]\times \{0\}$.

If we start with any Borel probability measure $\rho$ on $\mathbb{T}^2$,
$\mathcal{F}^n_*\rho$ must have a limit point by compactness. Define $\mu=(\mathcal{F}_*\rho)|_{(\mathbb{T}^2 \setminus [-a,a]\times \{0\})}$ and the initial
measure for the Markov chain $\tau$ such that $\tau[i]=\mu(L_i)$. Lemma \ref{lemma:mc vs dynamics modified} and Markov chain transience imply that all the limit points of $\mathcal{F}^{Nn}_*\mu$ are singular measures supported on $[-a,a] \times \{0\}$. Same conclusion can be applied if we define $\mu=(\mathcal{F}^k_*\rho)|_{(\mathbb{T}^2 \setminus [-a,a]\times \{0\})}$ for any $1 \geq k \geq N-1$. Therefore all the limit points of $\mathcal{F}^{n}_*\rho$ are singular measures supported on $[-a,a] \times \{0\}$. \ep \medbreak

\acks I would like to thank my Ph.D. thesis advisor Lai-Sang Young for fruitful discussions, effective criticism, and useful comments on many drafts of this paper.

\end{document}